\documentclass[11pt]{amsart}
\usepackage{amssymb,latexsym,amsmath,amsfonts,color}
\newtheorem{proposition}{Proposition}[section]
\newtheorem{theorem}[proposition]{Theorem}
\newtheorem{corollary}[proposition]{Corollary}
\newtheorem{example}[proposition]{Example}

\newtheorem{remark}[proposition]{Remark}

\numberwithin{figure}{section}

\newcommand{\bdpm}{\begin{displaymath}}
\newcommand{\edpm}{\end{displaymath}}
\newcommand{\beas}{\begin{eqnarray*}}
\newcommand{\eeas}{\end{eqnarray*}}
\newcommand{\ba}{\begin{array}}
\newcommand{\ea}{\end{array}}

\newcommand{\brm}{\begin{rm}}
\newcommand{\erm}{\end{rm}}

\newcommand{\s}[1]{\bar{#1}}
\newcommand{\Isim}{\stackrel{I}{\sim}}
\newcommand{\ssim}{\stackrel{s}{\sim}}

\title{Wilf classification of three and four letter\\ signed patterns}

\author{W. M. B. Dukes, T. Mansour and A. Reifegerste}

\address{School of Mathematical Sciences, University College Dublin, Ireland}
\email{dukes@maths.ucd.ie}
\address{Department of Mathematics, Haifa University, Israel}
\email{toufik@math.haifa.ac.il}
\address{Fakult\"{a}t f\"{u}r Mathematik und Physik, Universit\"at Hannover, Germany}
\email{reifegerste@math.uni-hannover.de}

\subjclass[2000]{Primary: 05A15; Secondary: 05A05}

\keywords{forbidden subsequences, pattern-avoiding permutations, signed permutations, signed involutions, Wilf equivalence}

\begin{document}

\begin{abstract}
We give some new Wilf equivalences for signed patterns which allow the complete classification of signed patterns of lengths three and four. The problem is considered for pattern avoidance by general as well as involutive signed permutations.
\end{abstract}
\maketitle

% INTRODUCTION

\section{Introduction}

Let $S_n$ and $B_n$ denote the symmetric group and hyperoctahedral group on $[n]=\{1,2,\ldots,n\}$, respectively. We regard the elements of $B_n$ as signed permutations written as a word $\pi=\pi_1\pi_2\ldots\pi_n$ in which each of the integers $1,2,\ldots,n$ appears, possibly signed which is represented by a bar. Throughout this paper, we use the terms ``barred/unbarred'' and ``negative/positive'' synonymously. Furthermore let $I_n$ and $SI_n$ be the set of involutions in $S_n$ and $B_n$, respectively.\\[2ex]
The {\it barring operation} maps $\pi\in B_n$ to the signed permutation $\s{\pi}$ which is obtained from $\pi$ by changing the sign of all elements. Clearly, we have $\s{\s{\pi}}=\pi$. The absolute value notation means $|\pi_i|=\pi_i$ if $\pi_i$ is not barred and $|\pi_i|=\s{\pi}_i$ otherwise. In particular, we write $|\pi|$ for the permutation $|\pi_1||\pi_2|\ldots|\pi_n|$.\\[2ex]
First we recall the concept of pattern-avoidance for signed permutations. A signed permutation $\pi\in B_n$ is said to {\it contain the pattern} $\tau\in B_k$ if there exists a sequence $1\le i_1<i_2<\ldots<i_k\le n$ such that $|\pi_{i_a}|<|\pi_{i_b}|$ if and only if $|\tau_a|<|\tau_b|$ and $\pi_{i_a}>0$ if and only if $\tau_a>0$ for all $1\le a,b\le k$. Otherwise, $\pi$ is called a {\it $\tau$-avoiding} permutation. By $M(\tau)$ we denote the set of all elements of $M$ which avoid the pattern $\tau$, and by $|M(\tau)|$ its cardinality.\\[2ex]
Two signed patterns $\sigma$ and $\tau$ are called Wilf equivalent, in symbols $\sigma\sim\tau$, if they are avoided by the same number of signed $n$-permutations for each $n\ge 1$. When we consider avoidance only in the set of signed involutions then we use $\Isim$ to denote the equivalence.\\
Obviously, if $\sigma,\tau\in S_k$ are Wilf equivalent patterns in $S_k$, i. e. $|S_n(\sigma)|=|S_n(\tau)|$ for all $n\ge1$,
then they are also equally restrictive for $B_n$.\\[2ex]
The classification of patterns up to Wilf equivalence is a basic problem in the theory of forbidden subsequences. Whereas many equivalences can be simply put down to symmetries within the symmetric and hyperoctahedral groups, detecting and proving {\it all} Wilf equivalences in $S_k$ or $B_k$ is a difficult task which requires a more subtle approach.
Over a period of 20 years, the classification of permutation patterns up to length 7 has been successfully completed.
The bijection between $S_n(123)$ and $S_n(132)$ given by Simion and Schmidt \cite{simion-schmidt} initiated this; it showed the whole of $S_3$ to be one Wilf class. The classification of $S_4$ turned out to be considerably more complicated and was done by West and Stankova in the series of papers \cite{west}, \cite{stankova1}, \cite{stankova2}. Using a stronger equivalence relation, Babson and West \cite{babson-west} showed that the patterns $(k-1)k\tau$ and $k(k-1)\tau$ for any $\tau\in S_{k-2}$ and $(k
 -2)(k-1)k\tau$ and $k(k-1)(k-2)\tau$ for any $\tau\in S_{k-3}$ are equally restrictive which provides the classification of $S_5$. However, the generalization $(k-l+1)(k-l+2)\ldots k\tau\sim k(k-1)\ldots(k-l+1)\tau$ for any $\tau\in S_{k-l}$ given by Backelin, West, and Xin \cite{backelin-west-xin} was not sufficient to classify $S_6$ completely. The last open case was settled by the equivalence of $(k-1)(k-2)k\tau$ and $(k-2)k(k-1)\tau$, shown in Stankova and West \cite{stankova-west}. With that it was even possible to complete the classification of patterns of length 7 which is the current limit.\\
Analogous investigations for pattern avoidance by involutive permutations were conducted in Guibert \cite{guibert}, Jaggard \cite{jaggard} and Bousquet-M{\'e}lou and Steingr{\'i}msson \cite{bousquet-steingrimsson}.\\[2ex]
In this paper the classification problem is completed for signed patterns of length at most 4.
Simion \cite{simion} proved that all elements of $B_2$ are Wilf equivalent to each other if pattern avoidance by all signed permutations is considered. For signed involutions, $B_2$ falls into two Wilf classes, see \cite{dukes-mansour}. We determine the Wilf classes for three and four letter signed patterns both regarding avoidance by general and involutive signed permutations.

% WILF EQUIVALENCES FOR SIGNED PATTERNS

\section{Wilf equivalences for signed patterns}

Before enumerating signed pattern-avoiding permutations, we shall study the equivalences between signed patterns. This will decrease considerably the number of patterns that need individual attention. There are many trivial equivalences based on symmetries. For a signed permutation $\pi=\pi_1\pi_2\ldots\pi_n$, we define its {\it reverse} as the permutation $\pi_n\pi_{n-1}\ldots\pi_1$ and its {\it complement} as the permutation whose $i$th element is $n+1-\pi_i$ if $\pi_i$ is positive and $-(n+1)-\pi_i$ otherwise. It is obvious that two signed patterns are Wilf equivalent if one of them can be transformed into the other by a sequence of reverse, complement or barring operations. (Note that the symmetry group is smaller if we consider the problem in the set of signed involutions. The group is then generated by the barring operation and the composition of reverse and complement operations.) Furthermore, every pattern is obviously Wilf equivalent to its inverse since $\tau$ is avoided by $\pi$ if and only if $\tau^{-1}$ is avoided by $\pi^{-1}$.\\[2ex]
In this section we will give some nontrivial equivalences which are the key to the complete classification of short signed patterns. Our bijective proofs use essentially the properties of right-to-left maxima of permutations. The idea to consider these special elements is based on the classical bijection given by \cite{simion-schmidt} for patterns of length three.\\
An element $\pi_i$ of a word $\pi$ is called a {\it right-to-left maximum} if it is greater than all elements that follow it, i.e. $\pi_i>\pi_j$ for all $j>i$. We define successively the $r$-right-to-left maxima for a signed permutation $\pi\in B_n$. Let $\pi^{(1)}$ be the subword consisting of all unbarred elements of $\pi$. For $r\ge 1$, the right-to-left maxima of $\pi^{(r)}$ are called {\it $r$-right-to-left maxima} of $\pi$ where $\pi^{(r+1)}$ is the subword obtained from $\pi^{(r)}$ by removing all $r$-right-to-left maxima.\\
For example, the signed permutation $\pi=2\:\s{5}\:6\:3\:10\:\s{8}\:4\:1\:7\:\s{9}\in B_{10}$ has the 1-right-to-left maxima 10 and 7; the 2-right-to-left maxima 6, 4, and 1; the 3-right-to-left maximum 3; and
the 4-right-to-left maximum 2.\\
Note that the $r$-right-to-left maxima of $\pi$ form a decreasing subsequence for each $r$. Furthermore, each $r$-right-to-left maximum $a$ is the initial term of an increasing subsequence of unbarred elements in $\pi$ of length $r$, and there is no increasing subsequence of unbarred elements of length $r+1$ which starts with $a$.\\[2ex]
Our first result states that barring the prefix of an increasing pattern yields an equivalence.

\begin{theorem} \label{prefix1}
For any pattern $\tau\in B_l$ and $k\ge l$, we have
\bdpm
\tau(l+1)(l+2)\ldots k\;\sim\;\s{\tau}(l+1)(l+2)\ldots k.
\edpm
Moreover, the patterns are also equivalent if we consider their avoidance by signed involutions.
\end{theorem}

\begin{proof}
We construct a bijection $\Phi_{k-l}$ from the set of all $\tau(l+1)(l+2)\ldots k$-avoiding permutations in $B_n$ to the set of all $\bar{\tau}(l+1)(l+2)\ldots k$-avoiding permutations in $B_n$ which preserves all $(k-l)$-right-to-left maxima.
In addition, $\Phi_{k-l}$ maps involutions to involutions again.\\
Given a signed permutation $\pi\in B_n(\tau(l+1)(l+2)\ldots k)$, we define $\sigma=\Phi_{k-l}(\pi)$ to be the permutation obtained from $\pi$ by barring all elements $\pi_i$ having a $(k-l)$-right-to-left maximum to their right which is greater than $|\pi_i|$. (This applies for an unbarred element $\pi_i$ if and only if $\pi_i$ is an $r$-right-to-left maximum with $r>k-l$.)\\
Why does the map leave all $(k-l)$-right-to-left maxima fixed? Let $a$ be any $(k-l)$-right-to-left maximum of $\pi$.
Then all elements following $a$ whose absolute value is greater than $a$ are left unchanged under $\Phi_{k-l}$ because the $(k-l)$-right-to-left maxima decrease. Thus, by construction, $a$ is also a $(k-l)$-right-to-left maximum of $\sigma$.
(By similar reasoning, it can be shown that all the $r$-right-to-left maxima for $r\le k-l$ are preserved.)\\
Since we have $|\pi_i|=|\sigma_i|$ for all $i$ and all the $(k-l)$-right-to-left maxima are preserved, $\Phi_{k-l}$ is an involution and hence bijective.\\
Why does $\sigma$ avoid the pattern $\bar{\tau}(l+1)(l+2)\ldots k$? We consider the $(k-l)$-right-to-left maximum $a$ of $\pi$ again. Note that an increasing unbarred subsequence of length $k-l$ starts in $a$. Consequently, there is no subsequence $\pi_{i_1},\ldots,\pi_{i_l}$ of elements preceding $a$ with $|\pi_{i_j}|<a$ for all $j$ which is order-isomorphic to $\tau$.
By applying $\Phi_{k-l}$, the element $a$ is unchanged while all elements with an absolute value smaller than $a$ to the left of $a$ are barred. Therefore there is no subsequence $\sigma_{i_1},\ldots,\sigma_{i_l}$ formed of these elements which is order-isomorphic to $\s{\tau}$.\\
Why is $\sigma$ an involution whenever $\pi$ is it? By construction, $\Phi_{k-l}$ has an effect on the sign but not on the position of the elements. Let $\pi_i=j$. (Clearly, $\pi_i$ and $\pi_j$ have the same sign since $\pi$ is an involution.)
If there is a $(k-l)$-right-to-left maximum $a$ to the right of $\pi_i$ with $a>|\pi_i|$ then the element $b=\pi^{-1}_a$ is a $(k-l)$-right-to-left maximum to the right of $\pi_j$ satisfying $b>|\pi_j|$. Thus $\Phi_{k-l}$ changes the sign of $\pi_i$ if and only if it does so for $\pi_j$.
\end{proof}

\begin{example}
\brm
Let $\pi=2\:\s{5}\:6\:3\:10\:\s{8}\:4\:1\:7\:\s{9}\in B_{10}$ again and $k-l=2$.
The following figure shows the effect of $\Phi_2$. We use the usual array representation of permutations with an
additional colouring to make a distinction regarding the element sign. Barred elements are represented by black points while unbarred elements are represented by white points. The $2$-right-to-left maxima are bordered. All the elements (points) which have to change their sign (colour) are contained in the grey region (the union of the south-west regions of all 2-right-to-left maxima).
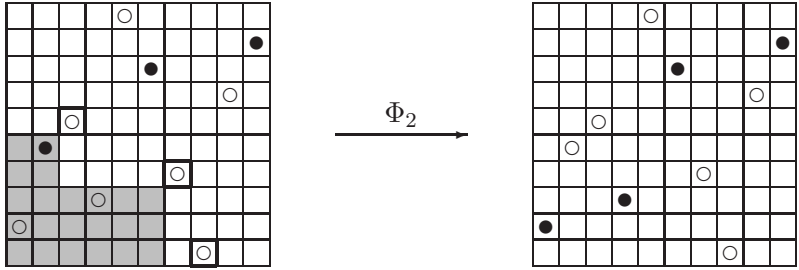
\begin{center}
\unitlength=0.35cm
\definecolor{g}{gray}{0.75}
\fboxsep0cm
\fboxrule0cm
\begin{picture}(30,10)
\put(0,0){\fcolorbox{g}{g}{\makebox(2,5){}}}
\put(2,0){\fcolorbox{g}{g}{\makebox(4,3){}}}
\multiput(0,0)(0,1){11}{\line(1,0){10}}
\multiput(0,0)(1,0){11}{\line(0,1){10}}
\put(0.5,1.5){\circle{0.5}}
\put(1.5,4.5){\circle*{0.5}}
\put(2.5,5.5){\circle{0.5}}
\put(3.5,2.5){\circle{0.5}}
\put(4.5,9.5){\circle{0.5}}
\put(5.5,7.5){\circle*{0.5}}
\put(6.5,3.5){\circle{0.5}}
\put(7.5,0.5){\circle{0.5}}
\put(8.5,6.5){\circle{0.5}}
\put(9.5,8.5){\circle*{0.5}}
\multiput(20,0)(0,1){11}{\line(1,0){10}}
\multiput(20,0)(1,0){11}{\line(0,1){10}}
\put(20.5,1.5){\circle*{0.5}}
\put(21.5,4.5){\circle{0.5}}
\put(22.5,5.5){\circle{0.5}}
\put(23.5,2.5){\circle*{0.5}}
\put(24.5,9.5){\circle{0.5}}
\put(25.5,7.5){\circle*{0.5}}
\put(26.5,3.5){\circle{0.5}}
\put(27.5,0.5){\circle{0.5}}
\put(28.5,6.5){\circle{0.5}}
\put(29.5,8.5){\circle*{0.5}}
\put(12.5,5){\vector(1,0){5}}
\put(15,5.75){\makebox(0,0)[cc]{$\Phi_2$}}
\linethickness{1pt}
\multiput(7,0)(0,1){2}{\line(1,0){1}}\multiput(7,0)(1,0){2}{\line(0,1){1}}
\multiput(6,3)(0,1){2}{\line(1,0){1}}\multiput(6,3)(1,0){2}{\line(0,1){1}}
\multiput(2,5)(0,1){2}{\line(1,0){1}}\multiput(2,5)(1,0){2}{\line(0,1){1}}
\end{picture}
\end{center}
We obtain $\sigma=\s{2}\:5\:6\:\s{3}\:10\:\s{8}\:4\:1\:7\:\s{9}\in B_{10}$ which has the same 2-right-to-left maxima as $\pi$.
\erm
\end{example}

\begin{corollary} \label{increasing pattern}
Let $\tau\in B_k$ with $|\tau|=12\ldots k$. Then $\tau\sim 12\ldots k$ and even $\tau\Isim 12\ldots k$.
\end{corollary}

\begin{proof}
Let $l$ be the maximal integer with $\tau_l=\s{l}$. We may assume that $l<k$; otherwise we consider $\s{\tau}$ which
is trivially Wilf equivalent to $\tau$. By Theorem \ref{prefix1}, the patterns $\tau_1\ldots\tau_l(l+1)(l+2)\ldots k$ and
$\s{\tau_1}\ldots\s{\tau_l}(l+1)(l+2)\ldots k$ are Wilf equivalent (also when considering signed involutions only).
Induction on $l$ yields the assertion.
\end{proof}

Our next result shows that barring the prefix of a pattern also yields an equivalent pattern if the monotone part is a decreasing sequence.

\begin{theorem} \label{prefix2}
For any pattern $\tau\in B_l$ and $k\ge l$, we have
\bdpm
\tau k(k-1)\ldots(l+1)\;\sim\;\s{\tau}k(k-1)\ldots(l+1).
\edpm
This relation remains true when we consider the pattern-avoiding signed involutions only.
\end{theorem}

\begin{proof}
This proof closely follows that of Theorem \ref{prefix1}, however the way in which we identify those elements whose sign will change differs slightly. Let $\pi\in B_n(\tau k(k-1)\ldots(l+1))$ and define $\sigma=\Psi_{k-l}(\pi)$ to be the permutation obtained from $\pi$ by changing the sign of all elements $\pi_i$ having (at least) $k-l$ 1-right-to-left maxima to their right which are greater than $|\pi_i|$. (In the graph of $\pi$, these are just the points $(i,\pi_i)$ for which there is a decreasing sequence of $k-l$ white points in the region north-east of $(i,\pi_i)$.)\\
By definition, there are no unbarred elements in the north-east of any 1-right-to-left maximum; thus these maxima are fixed under $\Psi_{k-l}$ (and hence also 1-right-to-left maxima of $\sigma$). Consequently, $\Psi_{k-l}(\sigma) = \pi$ and hence $\Psi_{k-l}$ is a bijection between $B_n(\tau k(k-1)\ldots(l+1))$ and $B_n(\s{\tau} k(k-1)\ldots(l+1))$.\\
Obviously, $\sigma$ avoids $\s{\tau}k(k-1)\ldots(l+1)$. If we suppose otherwise, then there exist a sequence
$a_1a_2\ldots a_k$ in $\sigma$ that is order-isomorphic to $\s{\tau}k(k-1)\ldots(l+1)$. (We may assume that $a_{l+1},\ldots,a_k$ are 1-right-to-left maxima of $\sigma$.) Since $a_{l+1}>a_{l+2}>\ldots>a_k>|a_i|$
for all $i\in[l]$, we must have $\s{a}_1\ldots\s{a}_la_{l+1}\ldots a_k$ in $\pi$, and this is an occurrence of $\tau k(k-1)\ldots(l+1)$, which is false. The proof remains true when $\sim$ is replaced by $\Isim$ since $\Psi_{k-l}$
maps signed involutions to signed involutions.
\end{proof}

\begin{example}
\brm
Consider the signed permutation
$\pi=8\:3\:\s{5}\:10\:2\:\s{1}\:7\:6\:\s{9}\:4\in B_{10}$ which avoids
the pattern $\s{1}\:2\:5\:4\:3$. Its 1-right-to-left maxima are 10, 7, 6, and 4.
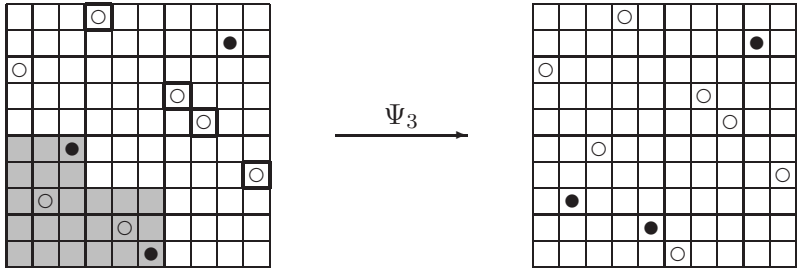
\begin{center}
\unitlength=0.35cm
\definecolor{g}{gray}{0.75}
\fboxsep0cm
\fboxrule0cm
\begin{picture}(30,10)
\put(0,0){\fcolorbox{g}{g}{\makebox(3,5){}}}
\put(3,0){\fcolorbox{g}{g}{\makebox(3,3){}}}
\multiput(0,0)(0,1){11}{\line(1,0){10}}
\multiput(0,0)(1,0){11}{\line(0,1){10}}
\put(0.5,7.5){\circle{0.5}}
\put(1.5,2.5){\circle{0.5}}
\put(2.5,4.5){\circle*{0.5}}
\put(3.5,9.5){\circle{0.5}}
\put(4.5,1.5){\circle{0.5}}
\put(5.5,0.5){\circle*{0.5}}
\put(6.5,6.5){\circle{0.5}}
\put(7.5,5.5){\circle{0.5}}
\put(8.5,8.5){\circle*{0.5}}
\put(9.5,3.5){\circle{0.5}}
\multiput(20,0)(0,1){11}{\line(1,0){10}}
\multiput(20,0)(1,0){11}{\line(0,1){10}}
\put(20.5,7.5){\circle{0.5}}
\put(21.5,2.5){\circle*{0.5}}
\put(22.5,4.5){\circle{0.5}}
\put(23.5,9.5){\circle{0.5}}
\put(24.5,1.5){\circle*{0.5}}
\put(25.5,0.5){\circle{0.5}}
\put(26.5,6.5){\circle{0.5}}
\put(27.5,5.5){\circle{0.5}}
\put(28.5,8.5){\circle*{0.5}}
\put(29.5,3.5){\circle{0.5}}
\put(12.5,5){\vector(1,0){5}}
\put(15,5.75){\makebox(0,0)[cc]{$\Psi_3$}}
\linethickness{1pt}
\multiput(3,9)(0,1){2}{\line(1,0){1}}\multiput(3,9)(1,0){2}{\line(0,1){1}}
\multiput(6,6)(0,1){2}{\line(1,0){1}}\multiput(6,6)(1,0){2}{\line(0,1){1}}
\multiput(7,5)(0,1){2}{\line(1,0){1}}\multiput(7,5)(1,0){2}{\line(0,1){1}}
\multiput(9,3)(0,1){2}{\line(1,0){1}}\multiput(9,3)(1,0){2}{\line(0,1){1}}
\end{picture}
\end{center}
The shaded region highlights the positions which are subject to sign-change. We obtain $\Psi_3(\pi)=8\:\s{3}\:5\:10\:\s{2}\:1\:7\:6\:\s{9}\:4\in B_{10}(1\:\s{2}\:5\:4\:3)$.
\erm
\end{example}

The third equivalence which we will use for the Wilf classification in the next sections is an immediate adjustment of a result for pattern of the symmetric group given by Backelin, West, and Xin \cite{backelin-west-xin}. They proved that the monotone suffix of any pattern can be reversed under Wilf equivalence. Bousquet-M{\'e}lou and Steingr{\'i}msson \cite{bousquet-steingrimsson} showed that this equivalence is preserved when considering avoidance by involutions.

\begin{theorem} \label{suffix reversal}
For any pattern $\tau\in B_l$ and $k\ge l$, we have
\bdpm
\tau(l+1)(l+2)\ldots k\;\sim\;\tau k(k-1)\ldots(l+1)
\edpm
where $\sim$ can be replaced by $\Isim$.
\end{theorem}

% CLASSIFICATION OF SIGNED PATTERNS OF LENGTH THREE

\section{Classification of signed patterns of length three}

There are 48 signed patterns of length three but by taking all symmetries in consideration, we can restrict our attention to the following six: $123$, $\s{1}23$, $1\s{2}3$, $132$, $\s{1}32$, $1\s{3}2$. It is well known
that all patterns of $S_3$ are avoided by the same number of permutations in $S_n$. Hence these patterns are also Wilf equivalent in $B_3$. This result and Theorem \ref{prefix1} and \ref{prefix2} reduce the number of patterns which
shall be considered up to two. (We have $123\sim \s{1}23\sim 1\s{2}3\sim 132\sim \s{1}32$.) A look at the initial terms of the sequences $(|B_n(\tau)|)_{n\ge 0}$
yields the Wilf classification of $B_3$.

\begin{table}[h]
\small
\begin{tabular}{|c|l|} \hline
123 & 1 2 8 47 358 3312 35784 440001\\\hline
$1\s{3}2$ & 1 2 8 47 358 3311 35738 438561\\\hline
\end{tabular}
\vspace*{1ex}
\caption{Wilf classes of $B_3$ (avoidance by $B_n$)}
\vspace*{-2ex}
\end{table}

Now we compare the cardinalities of the sets $SI_n(\tau)$ for $\tau\in B_3$. By symmetry and inversion, we obtain twelve patterns, namely $123$, $\s{1}23$, $1\s{2}3$, $132$, $\s{1}32$, $321$ (which are all equivalent by Theorem \ref{prefix1}, \ref{prefix2}, and  \ref{suffix reversal}), $1\s{3}2$, $231$, $\s{2}31$, $23\s{1}$, $\s{3}21$, $3\s{2}1$.

\begin{proposition}
We have $321\Isim 3\s{2}1$.
\end{proposition}

\begin{proof}
Let $\pi\in SI_n(321)$ and $\varphi$ be the map that changes the sign of all (necessarily barred) elements $\pi_i$ for which there are unbarred elements $\pi_j$ and $\pi_k$ with $j<i<k$ and $\pi_j>|\pi_i|>\pi_k$. (Therefore, $\pi_i$ is the middle element of an occurrence of $3\s{2}1$ in $\pi$.) Obviously, $\varphi(\pi)$ is a $3\s{2}1$-avoiding involution. Note that $jik$ is an occurrence of $3\s{2}1$ whenever $\pi_j\pi_i\pi_k$ is it. It is clear that $\varphi^2(\pi)=\pi$ and $\varphi$ is therefore a bijection.
\end{proof}

Using the same bijection, we can generalize the result as follows.

\begin{corollary} \label{prefix and suffix}
For any integers $k,s\ge1$ and any permutation $\tau$ on $\{s+1,s+2,\ldots,s+k\}$ we have
\bdpm
(2s+k)\ldots(s+k+1)\,\tau\,s\ldots1\Isim (2s+k)\ldots(s+k+1)\,\bar{\tau}\,s\ldots1.
\edpm
\end{corollary}

\begin{proof}
The case $s=1$ is proved by the bijection $\varphi$ established in the previous proof.
(Note that the elements $\pi_i$ need not be barred now.)
Clearly, $\pi\in SI_n$ avoids $(k+2)\,\tau 1$ if and only if $\varphi(\pi)$ avoids $(k+2)\,\bar{\tau} 1$. Induction on $s$ yields the assertion.
\end{proof}

\begin{remark}
\brm
The statement is also correct if we consider the relation $\sim$ instead of $\Isim$. (Of course, the map $\varphi$ can be applied to $B_n$ as well.) However, for general signed permutations this equivalence is already covered by Theorem \ref{prefix1} since $(k+2)\tau_2\ldots\tau_{k+1}1$ and $1\tau_{k+1}\ldots\tau_2(k+2)$ are related by symmetry. For instance, we have
\bdpm
53241 \ssim 14235 \stackrel{\ast}{\sim} \s{1}\s{4}\s{2}\s{3}5 \ssim \s{1}4235 \stackrel{\ast}{\sim} 1\s{4}\s{2}\s{3}5 \ssim
5\s{3}\s{2}\s{4}1
\edpm
where $\ssim$ stands for a symmetry relation and $\stackrel{\ast}{\sim}$ means equivalence by Theorem \ref{prefix1}.
\erm
\end{remark}

The computation of the initial terms of the sequences $(|SI_n(\tau)|)_{n\ge 0}$ for the remaining six patterns shows that there are no further equivalences.

\begin{table}
\small
\begin{tabular}{|c|l|} \hline
123 & 1 2 6 19 68 256 1032 4341 19154 87604 415868 \\\hline
$1\s{3}2$ & 1 2 6 20 74 288 1178 4978 21738 97420 448172 \\\hline
$231$ & 1 2 6 20 74 292 1220 5336 24316 114872 560840 \\\hline
$\s{2}31$ & 1 2 6 20 74 291 1207 5215 23362 107960 513236 \\\hline
$23\s{1}$ & 1 2 6 20 75 299 1259 5501 24813 114729 542074 \\\hline
$\s{3}21$ & 1 2 6 20 75 298 1250 5430 24347 111821 524921 \\\hline
\end{tabular}
\vspace*{1ex}
\caption{Wilf classes of $B_3$ (avoidance by $SI_n$)}
\vspace*{-2ex}
\end{table}

% CLASSIFICATION OF SIGNED PATTERNS OF LENGTH FOUR

\section{Classification of signed patterns of length four}

Extending the pattern length to four increases the complexity considerably. The 384 permutations in $B_4$ can be partitioned into 40 symmetry classes represented by:
\beas
&&1234,\;\s{1}234,\;1\s{2}34,\;\s{1}\s{2}34,\;\s{1}2\s{3}4,\;1\s{2}\s{3}4,\;1324,\;\s{1}324,\;1\s{3}24,\;\s{1}3\s{2}4,\;1\s{3}\s{2}4,\;2134,\;\s{2}134,\;21\s{3}4,\\
&&\s{2}\s{1}34,\;2\s{1}\s{3}4,\;\s{2}\s{1}\s{3}4,\;2143,\;\s{2}143,\;\s{2}\s{1}43,\;\s{2}1\s{4}3,\;2\s{1}\s{4}3,\;2314,\;\s{2}314,\;2\s{3}14,\;23\s{1}4,\;\s{2}\s{3}14,\;\s{2}3\s{1}4\\
&&2\s{3}\s{1}4,\;\s{2}\s{3}\s{1}4,\;2413,\;\s{2}413,\;\s{2}\s{4}13,\;2\s{4}\s{1}3,\;3214,\;\s{3}214,\;3\s{2}14,\;\s{3}2\s{1}4,\;3\s{2}\s{1}4,\;\s{3}\s{2}\s{1}4
\eeas
The study of the patterns in $S_4$ showed that they can be partitioned into only three Wilf classes, namely those of $1234$, $1324$, and $2314$. These equivalences are, of course, also valid in $B_4$. Applying Theorems \ref{prefix1}, \ref{prefix2}, and \ref{suffix reversal} reduces the number of patterns which we still need to consider to the following:
\beas
&&1234\sim\mbox{\scriptsize $3214\sim 2134\sim 2143\sim \s{1}234\sim 1\s{2}34\sim \s{1}\s{2}34\sim \s{1}2\s{3}4\sim 1\s{2}\s{3}4\sim \s{2}\s{1}34\sim 21\s{3}4 \sim \s{2}\s{1}\s{3}4\sim \s{2}\s{1}43\sim \s{3}\s{2}\s{1}4$},\\
&&1324\sim\mbox{\scriptsize $\s{1}324\sim 1\s{3}\s{2}4$\qquad(since $\s{1}324\sim \s{1}\s{3}\s{2}4$ by symmetry)},\\
&&1\s{3}24\sim\mbox{\scriptsize $\s{1}3\s{2}4$},\quad
\s{2}134\sim\mbox{\scriptsize $2\s{1}\s{3}4\sim \s{2}143$},\quad
2314\sim\mbox{\scriptsize $2413\sim \s{2}\s{3}\s{1}4$},\quad
\s{2}314\sim\mbox{\scriptsize $2\s{3}\s{1}4$},\quad
2\s{3}14\sim\mbox{\scriptsize $\s{2}3\s{1}4$},\\
&&23\s{1}4\sim\mbox{\scriptsize $\s{2}\s{3}14$},\quad
\s{3}214\sim\mbox{\scriptsize $3\s{2}\s{1}4$},\quad
3\s{2}14\sim\mbox{\scriptsize $\s{3}2\s{1}4$},\quad
\s{2}1\s{4}3,\quad
2\s{1}\s{4}3,\quad
\s{2}413,\quad
\s{2}\s{4}13,\quad
2\s{4}\s{1}3
\eeas

\begin{proposition}
We have $\s{2}1\s{4}3\sim 2\s{1}\s{4}3$.
\end{proposition}

\begin{proof}
Let $\pi\in B_n(\s{2}1\s{4}3)$. Define $\psi$ to be the map that changes the sign of all elements $\pi_i$ for which there is a sequence $\pi_j\pi_k$ to the right of $\pi_i$ with $|\pi_i|<|\pi_k|<|\pi_j|$ and $\pi_j<0<\pi_k$.
(That means, $\pi_i$ is the first element of an occurrence of $1\s{3}2$ or $\s{1}\s{3}2$ in $\pi$.) Obviously, $\psi(\pi)$ avoids $2\s{1}\s{4}3$ since the first two elements of any occurrence of $2\s{1}\s{4}3$ would be elements whose sign was changed by $\psi$ and hence $\pi$ would contain the pattern $\s{2}1\s{4}3$. Clearly, $\psi^2(\pi)=\pi$ and hence $\psi$ is bijective.
\end{proof}

\begin{table}[t]
\small
\begin{tabular}{|c|l||c|l|} \hline
$1234$ & 1 2 8 48 383 3798 44811 610354 &
$1324$ & 1 2 8 48 383 3798 44811 610355\\\hline
$1\s{3}24$ & 1 2 8 48 383 3798 44809 610214 &
$\s{2}134$ & 1 2 8 48 383 3798 44809 610280\\\hline
$\s{2}1\s{4}3$ & 1 2 8 48 383 3798 44810 610268 &
$2314$ & 1 2 8 48 383 3798 44810 610284\\\hline
$\s{2}314$ & 1 2 8 48 383 3798 44809 610212 &
$2\s{3}14$ & 1 2 8 48 383 3798 44809 610210\\\hline
$23\s{1}4$ & 1 2 8 48 383 3798 44809 610277 &
$\s{2}413$ & 1 2 8 48 383 3798 44809 610214\\\hline
$\s{2}\s{4}13$ & 1 2 8 48 383 3798 44808 610144 &
$2\s{4}\s{1}3$ & 1 2 8 48 383 3798 44808 610130\\\hline
$\s{3}214$ & 1 2 8 48 383 3798 44809 610279 &
$3\s{2}14$ & 1 2 8 48 383 3798 44809 610276\\\hline
\end{tabular}
\vspace*{1ex}
\caption{Wilf classes of $B_4$ (avoidance by $B_n$)}
\vspace*{-2ex}
\end{table}

Table 3 lists the initial terms of the sequences $(|B_n(\tau)|)_{n\ge0}$ for the fourteen patterns $\tau\in B_4$ which remain. For $n=7$, they are all different; hence the classification is done.\\[2ex]
Now we turn our attention to comparing the four letter patterns regarding their avoidance by signed involutions. Taking all symmetries for involutions into consideration, we obtain 78 classes. All the equivalences which we have obtained in the general case, apart from that one between $2314$ and $2413$ (given by Stankova \cite{stankova1}), are based on Theorem \ref{prefix1}, \ref{prefix2}, and \ref{suffix reversal}. Therefore they are still valid. By computation, we will see that
the equivalence of $2314$ and $2413$ gets lost for the involution case. So we may concentrate on the remaining cases:
\beas
&&1234,\;1324,\;1\s{3}24,\;\s{2}134,\;\s{2}1\s{4}3,\;2\s{1}\s{4}3,\;2314,\;\s{2}314,\;2\s{3}14,\;23\s{1}4,\;
2413,\;\s{2}413,\;\s{2}\s{4}13,\;2\s{4}\s{1}3,\\
&&\s{3}214,\;3\s{2}14;\\
&&2\s{4}13,\;\s{2}4\s{1}3,\;3412,\;\s{3}412,\;\s{3}\s{4}12,\;\s{3}4\s{1}2,\;3\s{4}\s{1}2,\;4123,\;\s{4}123,\;4\s{1}23,\;41\s{2}3,\;\s{4}\s{1}23,\;\s{4}1\s{2}3,\;4213,\\
&&\s{4}213,\;4\s{2}13,\;42\s{1}3,\;\s{4}\s{2}13,\;\s{4}2\s{1}3,\;4\s{2}\s{1}3,\;\s{4}\s{2}\s{1}3,\;4231,\;\s{4}231,\;4\s{2}31,\;\s{4}\s{2}31,\;4\s{2}\s{3}1,\;4312,\;\s{4}312,\\
&&4\s{3}12,\;43\s{1}2,\;\s{4}\s{3}12,\;\s{4}3\s{1}2,\;4321,\;\s{4}321,\;4\s{3}21,\;\s{4}\s{3}21,\;\s{4}3\s{2}1,\;4\s{3}\s{2}1.
\eeas
The patterns in the latter part represent classes arising when we only consider symmetries of involutions and have to been studied individually now.\\
Using a result of Guibert \cite{guibert}, we have $1234\Isim 3412\Isim 4321$. (Hence the unsigned patterns $\tau\in S_4$ contribute eight cases.) By Corollary \ref{prefix and suffix},
the relations $4321\Isim 4\s{3}\s{2}1$ and $4231\Isim 4\s{2}\s{3}1$ follow. Indeed, this completes the classification as
Table 4 shows. We obtain here fifty Wilf classes.\\
Note that the numbers $|SI_n(\s{2}1\s{4}3)|$ and $|SI_n(2\s{1}\s{4}3)|$ coincide for $n\le9$. This behaviour is quite exceptional; all the other patterns disclose their nonequivalence for shorter involutions.

\begin{table}[h]
\scriptsize
\begin{tabular}{|c|l||c|l|} \hline
$1234$ & 1 2 6 20 75 302 1299 5882 27899 137702 704716 &
$1324$ & 1 2 6 20 75 302 1299 5881 27889 137597 703878\\\hline
$1\s{3}24$ & 1 2 6 20 76 310 1354 6200 29644 146660 748752 &
$2\s{1}\s{4}3$ & 1 2 6 20 76 312 1378 6412 31246 157800 822452\\\hline
$\s{2}134$ & 1 2 6 20 76 310 1356 6224 29880 148592 763532 &
$\s{2}1\s{4}3$ & 1 2 6 20 76 312 1378 6412 31246 157800 822448\\\hline
$2314$ & 1 2 6 20 76 310 1358 6254 30202 151494 787398 &
$\s{2}314$ & 1 2 6 20 76 310 1358 6248 30117 150535 778460 \\\hline
$2\s{3}14$ & 1 2 6 20 76 310 1357 6238 30022 149808 773051 &
$23\s{1}4$ & 1 2 6 20 76 311 1368 6330 30676 154082 799383 \\\hline
$2413$ & 1 2 6 20 76 310 1360 6278 30444 153530 803578 &
$\s{2}413$ & 1 2 6 20 76 311 1370 6359 30994 156998 824015 \\\hline
$2\s{4}13$ & 1 2 6 20 76 311 1370 6358 30971 156682 820465 &
$\s{2}\s{4}13$ & 1 2 6 20 76 312 1381 6454 31678 161538 851968 \\\hline
$\s{2}4\s{1}3$ & 1 2 6 20 76 310 1359 6264 30290 152112 791459 &
$2\s{4}\s{1}3 $ & 1 2 6 20 76 312 1380 6442 31566 160672 845866 \\\hline
$\s{3}214$ & 1 2 6 20 76 311 1367 6318 30560 153147 792385 &
$3\s{2}14$ & 1 2 6 20 75 302 1300 5892 27993 138408 709859 \\\hline
$\s{3}412$ & 1 2 6 20 76 312 1378 6425 31428 159859 841636 &
$\s{3}\s{4}12$ & 1 2 6 20 76 312 1382 6476 31924 163898 871838 \\\hline
$\s{3}4\s{1}2$ & 1 2 6 20 75 302 1298 5868 27750 136364 693620 &
$3\s{4}\s{1}2$ & 1 2 6 20 76 312 1380 6452 31704 162232 860414\\\hline
$4123$ & 1 2 6 20 76 311 1368 6338 30797 155505 813216 &
$\s{4}123$ & 1 2 6 20 76 311 1370 6362 31015 157124 823967 \\\hline
$4\s{1}23$ & 1 2 6 20 76 311 1368 6337 30775 155205 809915 &
$41\s{2}3$ & 1 2 6 20 76 311 1369 6351 30924 156545 821054 \\\hline
$\s{4}\s{1}23$ & 1 2 6 20 76 311 1369 6350 30903 156262 817929 &
$\s{4}1\s{2}3$ & 1 2 6 20 76 311 1368 6336 30758 154992 807670 \\\hline
$4213$ & 1 2 6 20 76 310 1358 6252 30176 151212 784880 &
$\s{4}213$ & 1 2 6 20 76 311 1370 6359 30977 156715 820350 \\\hline
$4\s{2}13$ & 1 2 6 20 76 311 1368 6335 30754 155017 808670 &
$42\s{1}3$ & 1 2 6 20 76 312 1380 6444 31592 160973 848763 \\\hline
$\s{4}\s{2}13$ & 1 2 6 20 76 312 1380 6443 31573 160722 845999 &
$\s{4}2\s{1}3$ & 1 2 6 20 76 311 1368 6333 30719 154585 804274 \\\hline
$4\s{2}\s{1}3$ & 1 2 6 20 76 311 1369 6347 30866 155873 814600 &
$\s{4}\s{2}\s{1}3$ & 1 2 6 20 76 310 1358 6250 30140 150763 780284 \\\hline
$4231$ & 1 2 6 20 75 302 1299 5883 27911 137833 705870 &
$\s{4}231$ & 1 2 6 20 76 312 1379 6435 31510 160378 844431 \\\hline
$4\s{2}31$ & 1 2 6 20 75 302 1299 5882 27897 137674 704384 &
$\s{4}\s{2}31$ & 1 2 6 20 76 312 1378 6422 31380 159278 835774 \\\hline
$4312$ & 1 2 6 20 76 311 1368 6341 30840 155986 817676 &
$\s{4}312$ & 1 2 6 20 76 312 1380 6449 31661 161742 855816\\\hline
$4\s{3}12$ & 1 2 6 20 76 311 1369 6352 30936 156664 821993 &
$43\s{1}2$ & 1 2 6 20 76 311 1369 6351 30918 156433 819537 \\\hline
$\s{4}\s{3}12$ & 1 2 6 20 76 312 1382 6472 31872 163336 866840 &
$\s{4}3\s{1}2$ & 1 2 6 20 76 311 1368 6339 30804 155530 812915\\\hline
$\s{4}321$ & 1 2 6 20 76 311 1369 6347 30859 155752 813020 &
$4\s{3}21$ & 1 2 6 20 76 310 1358 6254 30200 151468 787094 \\\hline
$\s{4}\s{3}21$ & 1 2 6 20 76 312 1382 6468 31820 162774 861850 &
$\s{4}3\s{2}1$ & 1 2 6 20 76 310 1360 6274 30374 152658 794576\\\hline
\end{tabular}
\vspace*{1ex}
\caption{Wilf classes of $B_4$ (avoidance by $SI_n$)}
\vspace*{-2ex}
\end{table}
\newpage

% FINAL REMARKS

\section{Final remarks}

The application of all our results and the known Wilf classification of $S_5$ to suited representatives of the symmetry classes of signed patterns of length 5 yields 137 patterns whose (non)equivalence have to been proved. By computer checks up to $n=8$, it can be shown that these patterns form at least 58 different Wilf classes. It may be the case that the actual number of classes is much greater. By way of comparison, Table 5 lists the number of symmetry and Wilf classes in $B_k$ and $S_k$ for $k\le 5$, respectively:

\begin{table}[h]
\begin{tabular}{|l||c|c||c|c||c|c||c|c||c|c|} \hline
& $B_1$ & $S_1$ & $B_2$ & $S_2$ & $B_3$ & $S_3$ & $B_4$ & $S_4$ & $B_5$ & $S_5$\\\hline\hline
{\# symm. classes} & 1 & 1 & 2 & 1 & 6 & 2 & 40 & 7 & 284 & 23\\\hline
{\# Wilf classes} & 1 & 1 & 1 & 1 & 2 & 1 & 14 & 3 & $\ba{l}\ge 58\\[-1ex]\le 137\ea$ & 16\\\hline
\end{tabular}
\vspace*{1ex}
\caption{Number of symmetry and Wilf classes (avoidance by $B_n$ and $S_n$)}
\vspace*{-2ex}
\end{table}

Analogously, the five letter signed patters can be divided into at most 405 and at least 305 different classes when regarding their avoidance by signed involutions. The lower bound is obtained by comparing the size of $SI_9(\tau)$.

\begin{table}[h]
\begin{tabular}{|l||c|c||c|c||c|c||c|c||c|c|} \hline
& $B_1$ & $S_1$ & $B_2$ & $S_2$ & $B_3$ & $S_3$ & $B_4$ & $S_4$ & $B_5$ & $S_5$\\\hline\hline
{\# symm. classes} & 1 & 1 & 4 & 2 & 12 & 4 & 78 & 13 & 566 & 45\\\hline
{\# Wilf classes} & 1 & 1 & 2 & 1 & 6 & 2 & 50 & 8 & $\ba{l}\ge 305\\[-1ex]\le 405\ea$ & ?\\\hline
\end{tabular}
\vspace*{1ex}
\caption{Number of symmetry and Wilf classes (avoidance by $SI_n$ and $I_n$)}
\vspace*{-2ex}
\end{table}

% BIBLIOGRAPHY

\end{document}